%------------------------------------------------------------------------------
% divided2.tex  
%       Ringel-Zhang Pu: semi-Gorenstein-projective  Part III
% Beginn: 26.11.2019, revised Dec 6, 2019, revised March 29, 2020
%------------------------------------------------------------------------------
\magnification=\magstep1   
\input amstex
\UseAMSsymbols
\input pictex 
%\hoffset=0truecm \voffset=0truecm 
\vsize=23truecm
\NoBlackBoxes
\parindent=18pt
  
   \font\rmk=cmr8    \font\itk=cmti8

%\hrule height 2pt \vskip 3pt \hrule \bigskip\bigskip

\def\add{\operatorname{add}}
\def\Ker{\operatorname{Ker}}
\def\Cok{\operatorname{Cok}}
\def\soc{\operatorname{soc}}

\def\bdim{\operatorname{\bold{dim}}}

\def\top{\operatorname{top}}

  \def\ss{\ssize }
\def\arr#1#2{\arrow <1.5mm> [0.25,0.75] from #1 to #2}

\def\s{\hfill \square} 
%%%%%%%%%%%%%%%%%%%%%%%%%%%%%%%%%%%%%%%%%%%%%%%%%%%%%%%%%%
\vglue1truecm
\centerline{\bf  Koszul modules (and the $\Omega$-growth of modules)}
	\smallskip
\centerline{\bf over short local algebras.}
                     \bigskip
\centerline{Claus Michael Ringel, Pu Zhang}
                \bigskip\medskip

\noindent {\narrower Abstract: \rmk Following the well-established
terminology in commutative algebra, any (not necessarily commutative)
finite-dimensional local algebra $\ss A$ with radical $\ss J$
will be said to be {\itk short} provided $\ss J^3 = 0$. As in the commutative case,
also in general, the asymptotic behavior of the Betti numbers
of modules seems to be of interest. 
As we will see, there are only few possibilities for the growth of the
Betti numbers of modules. We generalize results which are known for
commutative algebras, but some of our results seem to be new also 
in the commutative case. 
	\medskip
\noindent 
2010 Math Subject classification. Primary 16G10, Secondary 13D07, 16E65, 16G50, 20G42.
	\medskip
\noindent 
Supported by NSFC 11971304.
\par}
 
	\bigskip
%%%%%%%%%%%%%%%%%%%%%%%%%%%%%%%%%%%%%%%%%%%%%%%%%%%%%%%%%%%%%%%%%%%
{\bf 1. Introduction.}
	\medskip
The modules to be considered 
are left modules of finite length over a finite-dimensional algebra $A$ (if not 
otherwise asserted). We denote by $|M|$ the length of the module $M$ 
and define
$t(M) = t_0(M)= |\top M|.$ For $n\in \Bbb N$, let $t_n(M) = t(\Omega^nM)$, where
$\Omega M = \Omega_A M$ is the first syzygy module of $M$ (as in commutative algebra [BH,L], 
one may call 
the numbers $t_n(M)$ the {\it Betti numbers} of $M$).
	\bigskip
{\bf 1.1. The $\Omega$-growth of a module.} We draw the attention 
to the asymptotic behavior of the Betti numbers $t_n(M)$ of a module $M$. 
If $M$ is a module, we consider the following numerical invariant
$$
 \gamma(M) =  \limsup_n \root n \of {t_n(M)}
$$
which we call the {\it $\Omega$-growth} of $M$. 
If $M$ has finite projective dimension, then $\gamma(M) = 0$;
otherwise $\gamma(M) \ge 1.$ 
Note that 
$$
 \gamma(M) = \limsup_n \root n \of {t_n(M)} = \limsup_n \root n \of {|\Omega^nM|}.
$$
This follows from the fact 
$$
 t_n(M) \le |\Omega^nM| \le |{}_AA|\cdot t_n(M)
$$
for all $n\ge 0$ (since, $t(N) \le |N|\le |{}_AA|\cdot t(N)$ for any module $N$).

If $A$ is a local algebra with simple module $S$, we define $\gamma_A = \gamma(S)$. 
	\medskip
{\bf Theorem 1.} {\it Let $A$ be a local algebra. Then $\gamma_A \le |{}_AJ|$
and $\gamma(M) = \gamma(\Omega M) \le \gamma_A$ for any module $M$. If $M$ is a module such that 
$S$ is a direct summand of $\Omega^n M$ for some $n \ge 0$,
then $\gamma(M) = \gamma_A.$}
	\bigskip
{\bf 1.2.} A local algebra $A$ with radical $J = J(A)$ 
is said to be {\it short} provided $J^3 = 0.$ 
Let $e = e(A) = |J/J^2|$ and, for $A$ being short, let
$a = a(A) = |J^2|$; then we call $(e(A),a(A))$ the {\it Hilbert-type} of $A$. 
As we have mentioned, the 
algebras to be considered will usually be local finite-dimensional $k$-algebras, where
$k$ is a field, with radical $J$. In addition, we usually will assume that $A$ is short
and that $A/J = k$, so that $S = k.$ 

An $A$-module has Loewy length at most 2 iff it is annihilated by $J^2$.
If $M$ is a module with Loewy length at most 2, we call
$\bdim M = (t(M),|JM|)$ (or its transpose, if we need to invoke matrix multiplication) 
the {\it dimension vector} of $M$.
Let us remark that $|M| = t(M) + |JM|$.
Recall from [RZ] that a module $M$ is said to be {\it bipartite} provided $\soc M = J M$.
A module has Loewy length at most $2$ if and only if it is the direct sum of
a bipartite and a semisimple module.
	\medskip 
{\bf 1.3. Koszul modules.} Let $A$ be a short local algebra of Hilbert type $(e,a).$
We have seen in [RZ] that the matrix 
$$
 \omega^e_a = \bmatrix e & -1 \cr
                    a & 0 \endbmatrix,
$$
controls (in some way) the change of the dimension vectors 
of modules of Loewy length at most 2, when we apply $\Omega = \Omega_A$. 
Always, the vectors $\bdim \Omega M$ and $\omega^e_a\bdim M$ are related by $\omega^e_a$, but
usually they differ slightly (see the Main Lemma of [RZ]; we recall it in 3.1).
A module $M$ of Loewy length at most 2 will be
said to be {\it aligned} provided $\bdim \Omega M = \omega^e_a\bdim M$. We study the
aligned modules very carefully in section 3. 

The main aim of the paper is to discuss the existence
and the structure of Koszul modules as defined by Herzog-Iyengar in [HI]. 
A short local
algebra is called a {\it left Koszul algebra} provided the simple module $S$ is
a Koszul module. Note that the projective modules are always Koszul modules. 
If $M$ is a Koszul module, 
then also $\Omega M$ is a Koszul module and has Loewy length at most 2 (since we assume,
as always, that $A$ is a short local algebra). 
In our setting, a module $M$ of Loewy length at most 2 is a {\it Koszul module} iff
$\bdim \Omega^n M = (\omega^e_a)^n\bdim M$ for all $n\ge 0$, thus iff all modules
$\Omega^n M$ are aligned, for $n\ge 0$ (see section 4).  
We usually will restrict the attention
to Koszul modules of Loewy length at most 2.

{\it If $M$ has Loewy length at most $2$ and 
$\Omega^n M$ is bipartite for all $n > 0,$ then $M$ is Koszul} (see 3.1). 
Thus, if $M$ is 
not Koszul, then $S$ is a direct summand of $\Omega^n M$ for some $n > 0$, therefore
Theorem 1 asserts that $\gamma(M) = \gamma_A.$  
The following theorem deals 
with the short local
algebras which have a non-projective Koszul module.
	\medskip 

{\bf Theorem 2.} {\it Let $A$ be a short local algebra of Hilbert type $(e,a).$
If there exists a non-zero Koszul module $M$ 
of Loewy length at most $2$, then the algebra is left Koszul, 
we have $a\le \frac14 e^2$ and $\gamma_A = \frac12(e+\sqrt{e^2-4a}\,).$
In addition, either 
$\gamma(M) = \gamma_A$ or else $a > 0$ and 
$\gamma(M) = \frac12(e-\sqrt{e^2-4a}\,).$} 
	\medskip
We recall that the {\it spectral radius 
$\rho(\omega)$} of a linear transformation $\omega\:\Bbb R^n \to \Bbb R^n$ is the 
maximum of the absolute values of the
(complex) eigenvalues of $\omega$. 
Note that for $a\le \frac14 e^2$, we have $\rho(\omega^e_a) = \frac12(e+\sqrt{e^2-4a}\,)$,
see 5.4. Thus, Theorem 2 asserts that the existence of a 
non-zero Koszul module $M$ implies that $\gamma_A = \rho(\omega^e_a).$
	\medskip
Theorem 2 provides a generalization of what Lescot [L] calls his key lemma:
the assertion that (for $A$ a commutative short local algebra with $\soc A = J^2$)
the existence of a non-zero Koszul module of Loewy length at most 2 implies that
$S$ is a Koszul module (see [L], 3.6).
	\bigskip
{\bf Theorem 3.} {\it Let $A$ be a short local algebra of Hilbert type $(e,a).$
Let $M$ be a non-zero module of Loewy length at most $2$ with $\gamma(M) < \gamma_A$. 
Then $M$ is a Koszul module, and the numbers 
$\gamma(M)$ and $\gamma_A$ are positive integers with
$$
 e = \gamma(M)+\gamma_A,\quad\text{and}\quad a = \gamma(M)\cdot\gamma_A.
$$
In particular, we have $0 < \gamma(M) < \frac12e < \gamma_A < e$ and
$e^2-4a = (\gamma_A-\gamma(M))^2$ (thus $e^2-4a$ is 
the square of a positive integer; in particular, positive). 
Also, $\bdim M$ is a multiple of $(1,\gamma_A)$ 
and $\bdim \Omega^n M = \gamma(M)^n\bdim M$ for all $n\in\Bbb N$.}
	\medskip

For example, let us look at the special case $e = 7.$ 
If there is a non-zero module $M$ 
of Loewy length at most $2$ with $\gamma(M) < \gamma_A,$ then $\gamma(M) = 1$ or $2$ or $3$, 
thus $a = 6,\ 10,\ 12,$ respectively (and $\bdim M$ is 
a multiple of $(1,6),\ (1,5),\ (1,4)$, respectively). 
Let us exhibit for $e = 7$ the graph of $\rho(\omega^e_a) = 
\frac12(e^2+\sqrt{e-4a}\,)$ as a function of $a$ (it contains the pairs $(6,6),\ (10,5),\
(12,4),$ they are marked by bullets $\bullet$), 
as well as (marked by small circles $\circ$) 
the three possible pairs $(a,\gamma(M))$, where $M$ is a Koszul module 
with $\gamma(M)  < \gamma_A,$ namely the pairs $(6,1),\ (10,2),\ (12,3)$:
$$
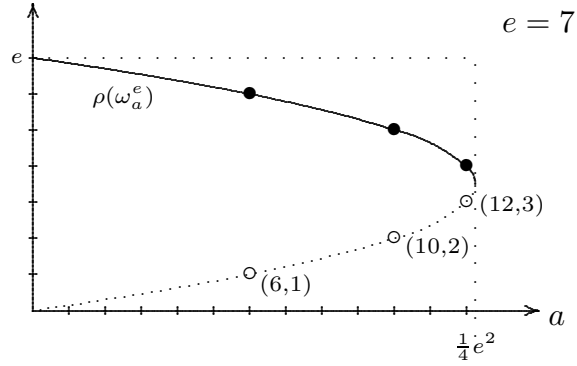
{\beginpicture
    \setcoordinatesystem units <.4cm,.4cm>
\put{$e = 7$} at 14 8
\multiput{} at 0 0  14 7.5 /
\put{$\ss \rho(\omega^e_a)$} at 2.5 5.9
\put{$\ss (6,1)$} at 7 0.7
\put{$\ss (10,2)$} at 11.2 1.75
\put{$\ss (12,3)$} at 13.25 2.9
\arr{0 0}{0 8.5}
\arr{0 0}{14 0}
\put{$a$} at 14.5 -.2
\plot 1 -.1  1 .1 /
\plot 2 -.1  2 .1 /
\plot 3 -.1  3 .1 /
\plot 4 -.1  4 .1 /
\plot 5 -.1  5 .1 /
\plot 6 -.1  6 .1 /
\plot 7 -.1  7 .1 /
\plot 8 -.1  8 .1 /
\plot 9 -.1  9 .1 /
\plot 10 -.1  10 .1 /
\plot 11 -.1  11 .1 /
\plot 12 -.1  12 .1 /
\plot -.1 1  .1 1 /
\plot -.1 2  .1 2 /
\plot -.1 3  .1 3 /
\plot -.1 4  .1 4 /
\plot -.1 5  .1 5 /
\plot -.1 6  .1 6 /
\plot -.1 7  .1 7 /
\put{$\ss\frac14 e^2$} at 12.25 -1
\put{$\ss e$} at -.45 7
\multiput{$\bullet$} at 6 6  10 5   12 4 /
\multiput{$\circ$} at 6 1   10 2   12 3   /
\setdots <2mm>
\plot 12.25 -.7 12.25  7  0 7 /
\setquadratic
\setdots <1mm>
\plot 0 0  6 1  10 2  11 2.35 12 3  12.2  3.25  12.25 3.5 /
\setsolid
\plot 0 7  6 6  10 5  11 4.65 12 4  12.2  3.75  12.25 3.5 /
\endpicture}
$$
%%%%%%%%%%%%%%%%%%%%%%%%%
	\bigskip
{\bf 1.4. Left Conca ideals.}
Let $A$ be a local algebra and $U$ an ideal of $A$. We say that $U$ is a
{\it left Conca} ideal provided $U^2 = 0$ and $J^2 \subseteq JU$.
If $A$ has a left Conca ideal $U$, then $A$ is short (namely, $J^3 \subseteq J^2U \subseteq 
JU^2 \subseteq U^2 = 0$). Since $J^2 \subseteq U$, the modules annihilated by $U$
have Loewy length at most 2. 
	\medskip 
{\bf Theorem 4.} 
{\it Let $A$ be a short local algebra. If $A$ has a left Conca ideal $U$,
then any module annihilated by $U$ is a Koszul module; in particular, $S$ is a Koszul module,
thus $A$ is a left Koszul algebra.}
	\medskip
This generalizes part of Theorem 1.1 of [AIS]. 
	\bigskip
{\bf 1.5. Construction of left Koszul algebras.}
	\medskip
{\bf Theorem 5.} {\it Given a pair $e,a$ of natural numbers, then the following assertions
are equivalent.
\item{\rm(i)} There is a short local algebra of Hilbert type (e,a) which is left Koszul.
\item{\rm(ii)} There is a commutative short local algebra of Hilbert type (e,a) which is
left Koszul.
\item{\rm(iii)} We have $a \le \frac14 e^2$.}
	\bigskip
{\bf Theorem 6.} {\it Let $c,d$ be positive integers. Let $e = c+d,\ a = cd.$ 
Then there are short local algebras of Hilbert type $(e,a)$ (even commutative ones)
with a Koszul module with dimension vector $(1,c).$}
	\medskip
Of course, if $c,d$ are positive integers and $e = c+d$ and $a = cd$, then we have
$a\le \frac14 e^2.$  
The algebras which we construct in the proof of Theorem 5 (showing that 
(iii) implies (ii)) and of Theorem 6 are short local algebras
with a left Conca ideal. 
	\bigskip
{\bf 1.6. A lower bound for $\gamma_A$.}
	\medskip
{\bf Theorem 7.} {\it Let $A$ be a short local algebra of Hilbert type $(e,a).$
If $a \le \frac14 e^2,$ then $\gamma_A \ge \frac12(e+\sqrt{e^2-4a}\,)$}.
	\medskip
In view of Theorems 5 and 2, the assertion of Theorem 7 can be strengthened
as follows. Let $\Cal A(a,e)$ be the class of all short local algebras of Hilbert
type $(e,a)$. Then: {\it For $a \le \frac14 e^2,$ the subset 
$\{\gamma_A\mid A\in \Cal A(e,a)\}$ of $\Bbb R$ 
has a minimal element, namely $\frac12(e+\sqrt{e^2-4a}\,).$}
(On the one hand, Theorem 7 shows that $\frac12(e+\sqrt{e^2-4a}\,)$ is a lower bound;
on the other hand, according to Theorem 5, there is a Koszul algebra $A$ in $\Cal A(e,a)$
and Theorem 2 asserts that  $\gamma_A = 
\frac12(e+\sqrt{e^2-4a}\,).$)

	\medskip
We have seen in [RZ] that there is a trichotomy
for short local algebras: There are the two special cases, first $a = 1$, 
second $a = e -1$, and then there are the remaining algebras with $a \notin \{1, e-1\}$ 
(for example, Gorenstein projective modules or non-zero minimal acyclic complexes
of projective modules do not exist if $a \notin \{1, e-1\}$). 
Theorem 5 yields a further separation: namely between $a \le \frac14e^2$ and
$a > \frac14e^2$: {\it The class $\Cal A(e,a)$ contains a Koszul algebra iff
$a \le \frac14 e^2$.} 
The disparity between $a \le \frac14e^2$ and
$a > \frac14e^2$
can be seen well if one looks at the spectral radius $\rho(\omega^e_a)$ as a function
of $a$ (fixing $e$): we have 
$\rho(\omega^e_a) = \frac12\bigl(e+\sqrt{e^2-4a}\,\bigr)$ for $a \le \frac14\,e^2$, 
and $\rho(\omega^e_a) = \sqrt a\ $ for $a \ge \frac14\,e^2$. 
$$
{\beginpicture
    \setcoordinatesystem units <1cm,1cm>
\arr{0 0}{5 0}
\arr{0 0}{0 2.5}
\plot -.1 1 .1 1 /
\plot -.1 2 .1 2 /
\plot 1 -.1 1 .1 /
\plot 2 -.1 2 .1 /
\plot 3 -.1 3 .1 /
\plot 4 -.1 4 .1 /
\put{$a$} at 5.2 -.1
\put{$\rho(\omega^e_a)$} at -.5 2.6
\put{$\ss e$} at -.3 2
\put{$\ss \frac12e$} at -.35 1
\put{$\ss \frac14e^2$} at 1 -.3
\put{$\ss e^2$} at 4 -.3
\setdots <2mm>
\plot 0 2  4 2  4 0 /
\plot 0 1  4 1 /
\plot 1 0  1 2 /
\setdots <.5mm>
\setquadratic
\plot 0 0  0.2 0.4  1 1  2.5  1.6  4 2 /
\setsolid
\plot 1 1  2.5  1.6  4 2 /
\plot 0 2  0.76  1.5  1 1 /
\endpicture}
$$
Note that if $a \le \frac14e^2$, thus Theorem 7 asserts that 
$\rho(\omega^e_a)$ is a lower bound for $\gamma_A$, and it seems that this is also 
true for $a > \frac14e^2$. 

	\bigskip
{\bf 1.7. Outline of the paper.} Sections 3 and 4 provide
characterizations of the aligned modules and the Koszul modules, respectively.
The $\Omega$-growth of modules is discussed in sections 2 and 5; in section 2,
there is the proof of Theorem 1, in section 5 the proof of Theorems 2 and 3.
Section 6 deals with left Conca ideals and presents the proof of Theorem 4.
In section 7 we construct suitable algebras with left Conca ideals 
in order to establish Theorems 5 and 6.
The final section 8 provides a lower bound for $\gamma_A$,
provided $a\le \frac14 e^2$. 

	\bigskip\bigskip
%%%%%%%%%%%%%%%%%%%%%%%%%%%%%%%%%%%%%%%%%%%%%%%%%%%%%%%%%%%%%%%%%%%%%%
%%%%%%%%%%%%%%%%%%%%%%%%%%%%%%%%%%%%%%%%%%%%%%%%%%%%%%%%%%%%%%%%%%%%%%
{\bf 2. The $\Omega$-growth of a module.}
	\medskip
First, let us consider an arbitrary finite-dimensional algebra $A$.
	\smallskip
{\bf 2.1. Lemma.} {\it If $M'$ is a direct summand of $M$, then $\gamma(M') \le \gamma(M).$}
	\medskip
Proof. If $M'$ is a direct summand of $M$, then $\Omega^n M'$ is a direct summand of $\Omega^n M,$
thus $|\Omega^nM'| \le |\Omega^nM|$ for all $n\ge 0.$ $\s$
	\medskip 
{\bf 2.2. Lemma.}
{\it If $0 \to M' \to M \to M'' \to 0$ is an exact sequence, then }
$$
 \gamma(M) \le \max\{\gamma(M'),\gamma(M'')\}.
$$
	\medskip
Proof. We start with minimal projective resolutions of $M'$ and $M''$. The horseshoe
lemma provides a (not necessarily minimal) projective resolution of $M$. 
This shows that $t_n(M) \le t_n(M')+t_n(M'')$ for all $n\ge 0.$ Therefore
$$
 \limsup \root n \of {t_n(M)} \le \max\{\limsup \root n \of {t_n(M')},
 \limsup \root n \of {t_n(M'')}\}.
$$
 $\s$
	\medskip
{\bf 2.3. Lemma.} {\it Let $A$ be a finite-dimensional local algebra and $M$ a module. Then
$\gamma(M) = \gamma(\Omega M) \le |{}_AJ|.$} 
	\medskip
Proof. Let $J$ be the radical of $A$ and $d = |{}_AJ|.$ Let $m$ be the Loewy length of ${}_AJ$
(thus $J^m \neq 0$, and $J^{m+1} = 0)$. Let $c = |{}_AJ^m|$.

For $n\ge 1$, the module $\Omega^nM$ is a submodule of $JP(\Omega^{n-1}M)$, thus it has 
Loewy length at most $m$. But if $N$ is a module of Loewy length at most $m$, then 
$J^mP(N) \subseteq \Omega N$. This shows that for $n\ge 1$, we have
$$
 J^mP(\Omega^nM) \subseteq \Omega^{n+1}M \subseteq JP(\Omega^nM).
$$
Since $|J^mP(\Omega^nM)| = ct_n(M)$ and $|JP(\Omega^nM)| = dt_n(M)$, we get
$$
 ct_n(M) \le |\Omega^{n+1}M| \le dt_n(M),
$$
thus 
$$
 \gamma(M) = \limsup \root n \of {t_n(M)} =
 \limsup \root n \of {|\Omega^{n+1}M|} = \gamma(\Omega M).
$$

On the other hand, we have by induction $t_{n+1}(M) \le d^nt_1(M)$, therefore
$$
 \gamma(\Omega M) =
 \limsup \root n \of {t_{n+1}(M)} \le
 \limsup \root n \of {d^n t_1(M)} = d = |{}_AJ|.
$$
$\s$
	\bigskip
{\bf 2.4. Proof of Theorem 1.} We have seen in 2.3 that $\gamma_A = \gamma(S) \le |{}_AJ|$
and that $\gamma(M) = \gamma(\Omega M).$

It follows from 2.2 that $\gamma(M) \le \gamma(S)$, using induction on the length of
$M$. Thus $\gamma(M) \le \gamma_A$. 

Now assume that $S$ is a direct summand of $\Omega^nM$ for some $n\ge 0.$ 
Using 2.1 and 2.3, we get
$\gamma(M) \le \gamma(S) \le \gamma(\Omega^nM) = \gamma(M).$  $\s$
	\bigskip
{\bf 2.5. Remark.} We should stress that the $\Omega$-growth $\gamma(M)$ of a module $M$
measures the {\bf exponential} growth of the Betti numbers. A similar, but deviating
measure, the complexity, was introduced by Alperin and Evens [AE] 
in 1981 dealing with representations of a finite
group $G$: The {\it complexity} of a $kG$-module is the least integer $c$ 
such that there is a constant $\kappa > 0$ with $t_n(M) \le \kappa\cdot n^{c-1}$ for all 
$n\ge 1$. 
In contrast to the $\Omega$-growth, the complexity measures the {\bf polynomial} growth of
the Betti numbers. There is the following obvious observation:
{\it If a $kG$-module $M$ has finite complexity and  
is not projective, then $\gamma(M) = 1.$} This follows from the fact that $\lim_n \root n \of n = 1.$ 

Dealing with an arbitrary finite-dimensional algebra, it may be advisable to look 
at various measures for the growth of the Betti numbers. However, the present
investigation seems to indicate that for short local algebras, it is the $\Omega$-growth
as defined in the introduction which is the decisive invariant. 
	\bigskip\bigskip 
%%%%%%%%%%%%%%%%%%%%%%%%%%%%%%%%%%%%%%%%%%%%%%%%%%%%%%%%%%%%%%%%%%%%%%
{\bf 3. Aligned modules.}
	\medskip
From now on, $A$ will be a short local $k$-algebra with radical $J$ such that $S = A/J = k.$
	\medskip 
{\bf 3.1.} We recall from [RZ] the Main Lemma.
{\it Let $A$ be a short local algebra of Hilbert type $(e,a)$.
If $M$ is a module of Loewy length at most $2$, then there is a natural number $w$
such that 
$$
 \bdim \Omega M = \omega^e_a \bdim M + (w,-w),
$$
and such that $\Omega M$ has a direct summand of the form $S^w$.} 
	\medskip
According to the Main Lemma, we have $\bdim \Omega M = \omega^e_a\bdim M$ 
provided $\Omega M$ is bipartite. But this formula is valid for a larger class
of modules, namely the aligned modules. We are going to provide several
equivalent conditions for a module to be aligned.
	\medskip
{\bf 3.2.}
If $M$ is a module of Loewy length at most $2$, let $p\:P(M) \to M$
be a projective cover. We consider $\Omega M$ as a submodule
of $JP(M)$ with inclusion map $u\:\Omega M \to JP(M)$ and obtain in this
way the exact sequence 
$$
 \eta_M = (\ 0 @>>> \Omega M @>u>> JP(M) @>p>> JM @>>> 0\ ).
$$
(In order to see that this sequence is exact, we apply the Snake Lemma to the
following commutative diagram with exact rows:
$$
\CD
 0 @>>> JP(M)   @>>>   P(M) @>>>    \top  M @>>> 0 \cr
 @.      @VVV           @VVpV              @VVp''V \cr
 0 @>>> JM      @>>>   M    @>>>     \top M @>>> 0 \cr
\endCD
$$
The kernel of $p\:P(M) \to M$ is $\Omega M$. Since $p''$ is an isomorphism, and
$p$ is surjective, we see that the cokernel of $JP(M) \to JM$ is zero.)

Considering the top of the modules, the exact sequence $\eta_M$ 
yields the exact sequence
$$
 \overline \eta_M = (\ \top\Omega M @>\overline u>> \top JP(M) 
 @>\overline p>> JM @>>> 0\ ),
$$ 
here we use that $JM$ is semisimple, since the Loewy length of $M$ is at most 2.
	\medskip
{\bf 3.3. Proposition.} {\it Let $M$ be a module of Loewy length at most $2$. The
following conditions are equivalent.
\item{\rm(i)} $M$ is aligned (by definition, this means that
 $\bdim \Omega M = \omega^e_a\bdim M$).
\item{\rm(ii)} $t(\Omega M) = et(M)-|JM|.$
\item{\rm(iii)} $|J\Omega M| = at(M).$
\item{\rm(iv)} $J\Omega M = J^2P(M).$
\item{\rm(v)} $J^2P(M) \subseteq J\Omega M.$
\item{\rm(vi)} The inclusion map $u$ yields an injective map 
  $\overline u\: \top \Omega M \to \top JP(M)$. 
\item{\rm(vii)} The sequence $\eta_M$ induces an exact sequence
 $0 \to \top \Omega M @>>> \top JP(M) @>>> JM \to 0.$
\item{\rm(viii)} A minimal projective
  presentation $P_1 \to P_0 \to M \to 0$ induces an exact sequence
 $0 \to \top P_1 @>>> \top JP_0 @>>> JM \to 0.$
\par}
 	\medskip
Proof. Let us start with the equivalence of (ii), (iii).
The Main Lemma (see 3.1) asserts that
$$
 (t(\Omega M),|J\Omega M|) = (et(M)-|JM|+w,at(M)-w)
$$
for some $w$. Thus, if $t(\Omega M) = et(M)-|JM|$ (the condition (ii)), then $w = 0$
and therefore $|J\Omega M| = at(M)$ (the condition (iii)). And conversely, if the
condition (iii) is satisfied, then again we have $w = 0$, thus condition (ii) is satisfied.

Assertion (i) is the conjunction of (ii) and (iii), thus it is of course
equivalent to (i) and to (iii). 
The inclusion map $J\Omega M \subseteq J^2P(M)$ shows  
that (iii) and (iv) are equivalent. Since $\Omega M \subseteq JP(M)$, we always have
$J\Omega M \subseteq J^2P(M)$. Thus (iv) and (v) are equivalent.

For the equivalence of (iv) and (vi), we apply the Snake Lemma to the
following commutative diagram with exact rows
$$
\CD
 0 @>>> J\Omega M   @>>> \Omega M @>>> \top \Omega M @>>> 0 \cr
 @.      @VVu'V           @VVuV              @V\overline uVV \cr
 0 @>>> J^2P(M) @>>> JP(M) @>\pi>>     \top JP(M) @>>> 0 \cr
\endCD
$$
We have $\Ker(u) = 0$ and $\Cok(u) = JM$. Also, 
the vertical map $\overline u$ on the right is part of the exact sequence 
$\overline \eta_M$, thus its cokernel is also $JM$. Altogether, the Snake Lemma yields
the exact sequence $0 \to \Ker(\overline u) \to \Cok u' \to JM \to JM \to 0.$
The surjective map $JM \to JM$ has to be an isomorphism, thus 
the nap  $\Ker(\overline u) \to \Cok u'$ has to be
an isomorphism. This mean that $u'$ is surjective (the condition (iv))
if and only if $\overline u$ is injective (the condition (vi)). 

The conditions (vi) and (vii) are of course equivalent, since 
$\eta_M$ induces the exact sequence 
$\overline \eta_M = (\ \top\Omega M @>\overline u>> 
\top JP(M) @>\overline p>> \top JM @>>> 0\ ), $ and this is a short exact
sequence if and only if $\overline u$ is injective. 

The assertions (vii) and (viii) are equivalent, since $P_1 = P(\Omega M)$ 
and in this way, $\top P_1$ is identified with $\top \Omega M.$ 
$\s$
	\medskip
{\bf 3.4. Remark.} {\it Let $V$ be a proper left ideal of $A$. Then
${}_AA/V$ is aligned if and only if $J^2 \subseteq JV.$} Namely, 
$\Omega({}_AA/V) = V$, thus we deal with condition (iv).

In particular: 
{\it The simple module $S$ is always aligned,} since here we have $V = J.$ 
(Of course, we also may look at $\bdim S = (1,0);$ we have $\Omega S = {}_AJ$
and $\bdim {}_AJ = (e,a) = \omega^e_a(1,0),$ this is condition (i).)
	\medskip
{\bf 3.5.} We recall from [RZ], 13.2: 
{\it Let $A$ be a short local algebra and 
$M$ a module of Loewy length at most $2$.  If $\Omega M$ is bipartite,
then $M$ is aligned. Conversely, if $J^2 = \soc {}_AA$, and $M$ is aligned,
then $M$ is bipartite.}
	\medskip
The condition $J^2 = \soc{}_AA$ has been discussed quite carefully in section 13 of [RZ].
	\medskip

{\bf 3.6. Examples of aligned modules $M$ such that $\Omega M$ is not bipartite.}
	\smallskip
{\bf (1)} 
Note that 3.4 provides such an example, namely $S$ is always aligned, whereas
$\Omega S = {}_AJ$ is bipartite iff $J^2 = \soc{}_AA$.
Note that for all short local algebras with $a = 0$ 
and $e\ge 1$, but also for many other short local algebras, we have $J^2 \neq \soc {}_AA$. 
		\medskip
Here are two additional examples of indecomposable 
modules $M$ of Loewy length $2$ which are aligned, but $\Omega M$ is not bipartite.
	\smallskip
{\bf (2)} Here is a typical example of a short local algebra with $J^2\neq \soc{}_AA$:
the algebra $A$ generated by $x,y,z$, with 
relations $x^2, y^2, z^2, xy-yx,
xz, zx, yz, zy.$ 

Let $M = Ay \simeq A/(A y+A z),$ thus $\Omega M = A y \oplus A z$.
Then $\bdim M = (1,1)$
and  $\bdim \Omega M = (2,1).$ 
$$
{\beginpicture
    \setcoordinatesystem units <.8cm,.8cm>
%%%%%%%%%%%%%%%%%%%%%%%%%
\put{\beginpicture
\put{$x$\strut} at 0 1
\put{$y$\strut} at 1 1
\put{$z$\strut} at 2 1
\put{$yx$\strut} at 0 0
\arr{0 0.8}{0 0.2}
\arr{0.8 0.8}{0.2 0.2}
\multiput{$\ss x$} at 0.7 .4   /
\multiput{$\ss y$} at  -.25 .5   /
\put{$J$} at -1.2 1 
\endpicture} at -1 0
%%%%%%%%%%%%%%%%%%%%%%%%%
\put{\beginpicture
\multiput{$\bullet$} at 0 0  1 1 /
\arr{0.9 0.9}{0.1 0.1}
\multiput{$\ss x$} at .3 .6   /
\put{$M$} at -.5 1 
\endpicture} at 4 0
%%%%%%%%%%%%%%%%%%%%%%%%%
\put{\beginpicture
\multiput{$\bullet$} at 0 0  1 1  2 1 /
\arr{0.9 0.9}{0.1 0.1}
\multiput{$\ss x$} at .3 .6   /
\put{$\Omega M$} at -.5 1
\endpicture} at 8 0
\endpicture}
$$
	\smallskip
{\bf (3)} Consider now the algebra $A$
$A$ generated by $x,y,z$ with the relations 
$yx-xy, zy-x^2, zx, y^2, xz, yz, z^2$ (thus $J^2$ has the basis $x^2, yx$). 
We define $M$ by taking a suitable submodule $U$ of a projective module $P$ 
and define $M = P/U$ so that $\Omega M = U$.
Namely, let $P = A^3$ and let $\Omega M$ be the submodule of $P$
generated by $(x,y,0),\ (0,x,y),\ (0,z,0)$.
$$
{\beginpicture
  \setcoordinatesystem units <1cm,1cm>
%%%%%%%%%%%%%%%%%%%%%%%%%
\put{\beginpicture
    \setcoordinatesystem units <1cm,1.3cm>
\put{$x$\strut} at 0 1
\put{$y$\strut} at 1 1
\put{$z$\strut} at 2 1
\put{$x^2$\strut} at 0 -.0
\put{$yx$\strut} at 1 -.0
\arr{0 0.8}{0 0.2}
\arr{0.8 0.8}{0.2 0.2}
\arr{1 0.8}{1 0.2}
\arr{0.2 0.8}{0.8 0.2}
\multiput{$\ss x$\strut} at -.2 0.55  1.2 .55  /
\multiput{$\ss y$\strut} at .4 0.85   /
\multiput{$\ss z$\strut} at .61 0.85  /
\put{$J$} at -1.1 1 
\endpicture} at 0 0
%%%%%%%%%%%%%%%%%%%%%%%%%
\put{\beginpicture
    \setcoordinatesystem units <.9cm,1cm>
\plot 0 0  1 1  1 0  2 1  2 0  3 1  3 0 /
\multiput{$\bullet$} at 0 0  1 1  1 0  2 1  2 0  3 1  3 0  -1.5 0  4.5 0 /
\multiput{$\ss x$} at  0.4 0.6  1.4 0.6  2.4 0.6 /
\multiput{$\ss y$} at  0.8 0.3  1.8 0.3  2.8 0.3 /
\multiput{$\ss z$} at  -0.75 0.7  4 0.7 /
\arr{-1.3 0.12}{-1.4 0.05}
\arr{0.1 0.1 }{0.05 0.05}
\arr{1 0.1}{1 0.05}
\arr{1.1 0.1 }{1.05 0.05}
\arr{2 0.1}{2 0.05}
\arr{2.1 0.1 }{2.05 0.05}
\arr{3 0.1}{3 0.05}
\arr{4.4 0.1}{4.5 0.04}
\setquadratic
\plot -1.5 0  0 0.8  1 1 /
\plot 3 1  3.5 0.8  4.5 0 /
\put{$M$} at -1.5 1
\endpicture} at 6.5 0
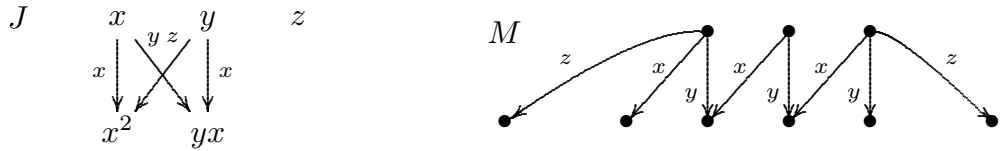
\endpicture}
$$
It turns out that $M$ is indecomposable and of 
Loewy length $2$. Since $A$ has
Hilbert type $(3,2)$ and $\bdim \Omega M = \bdim M = (3,6)$, 
we see that $M$ is aligned (the condition (i) is satisfied).
But $U = \Omega M$ is not bipartite: the submodule $J\Omega M$ is generated by 
$$
\gather 
x(x,y,0) = (x^2,xy,0),\
y(x,y,0) = (xy,0,0),\
z(x,y,0) = (0,x^2,0),\cr
x(0,x,y) =  (0,x^2,xy),\
y(0,x,y) = (0,xy,0),\
z(0,x,y) = (0,0,x^2),
\endgather
$$
thus equal to $J^2P$, and therefore $\Omega M$
is isomorphic to the direct sum of two copies of $A/J^2$ and one copy of 
the simple module $S$. Altogether, we see that $M$ is aligned, 
but $\Omega M$ is not bipartite. 
	\smallskip
Of course, in all the examples of 3.6, we have $J^2 \neq \soc{}_AA$, see 3.5. 
	\bigskip\bigskip
%%%%%%%%%%%%%%%%%%%%%%%%%%%%%%%%%%%%%%%%%%%%%%%%%%%%%%%%%%%%%%%%%%%%%%%%
{\bf 4. Koszul modules.}
	\medskip
{\bf 4.1. Koszul modules.} Following Herzog-Iyengar [HI] (see also [AIS]) 
a module $M$ will be said to be a {\it Koszul module} provided a minimal projective resolution 
$$
 \cdots \longrightarrow \quad P_{n+1} \quad \to \quad P_n \quad \to \quad P_{n-1} 
  \quad \to \cdots  \to 
    \quad P_0 \quad \to \quad M \quad \to \quad 0 \qquad\qquad
$$
induces for any $n\ge 0$ an exact sequence
$$
 \epsilon^M_n\:\qquad
  0 \to P_n/JP_n \to JP_{n-1}/J^2P_{n-1} \to \cdots \to 
    J^nP_0/J^{n+1}P_0 \to J^nM/J^{n+1}M \to 0
$$
(note that the image of $d_{i+1}\:P_{i+1} \to P_i$
is contained in $JP_i,$ thus $d_{i+1}(J^nP_{i+1}) \subseteq J^{n+1}P_i$ for all $n\ge 0$). 

A local algebra is called a {\it left Koszul algebra} provided the simple module $S$
is a Koszul module. 
	\medskip
Of course, the projective modules are always Koszul modules. If $M$ is a Koszul module, 
then also $\Omega M$ is a Koszul module and has Loewy length at most 2 (since we assume,
as always, that $A$ is a short local algebra). In the following, we usually will 
restrict the attention
to Koszul modules of Loewy length at most 2.
	\medskip
{\bf Proposition.} {\it Let $A$ be a short local algebra and $M$ a module
of Loewy length at most $2$.
The following conditions are equivalent:
\item{\rm(i)} $M$ is a Koszul module.
\item{\rm(ii)} For every $n\ge 1$, the exact sequence 
  $0 \to \Omega^{n}M \to P(\Omega^{n-1}M) \to \Omega^{n-1}M \to 0$ induces
  an exact sequence $0 \to \Omega^{n}M/J\Omega^{n}M \to 
  JP(\Omega^{n-1}M)/J^2P(\Omega^{n-1}M) \to J\Omega^{n-1}M \to 0.$
\item{\rm(iii)} $\bdim \Omega^nM = (\omega^e_a)^n\bdim M$ for all $n\ge 0.$
\item{\rm(iv)} The modules $\Omega^n M$ with $n\ge 0$ are aligned.}
	\medskip
Proof of the equivalence of (i) and (ii). We use the 
isomorphisms $\top P(\Omega^iM)\to \top \Omega^iM$.
Also note that $J^2\Omega^i M = 0$ for all $i\ge 0.$

We can rewrite $\epsilon^M_n$ as:
$$
 \epsilon^M_n\:\qquad
  0 \to \top P(\Omega^nM) \to \top JP(\Omega^{n-1}M) \to \cdots \to 
    \top J^nP(M) \to \top J^nM \to 0
$$
For $n=0$, this sequence $\epsilon^M_0\: 0 \to \top P(M) \to \top M \to 0$ 
is always exact. 

For $n\ge 1$, we can use the isomorphism $\top P(\Omega^nM) \to \top\Omega^nM$
in order to rewrite $\epsilon^M_n$ as
$$
 0 \to \top\Omega^{n}M \to 
  \top JP(\Omega^{n-1}M) \to J\Omega^{n-1}M \to 0.
$$
Note that this is just the exact sequence $\overline \eta_{\Omega^nM}$
as considered in 3.2.

The equivalence of (ii) and (iii) for $n\ge 1$ is given by Proposition 3.3,
namely we use the equivalence of (i) and (vii) for $M$ replaced by $\Omega^{n-1}M$.
By the definition of an aligned module, the conditions (iii) and (iv) are the same.
$\s$
	\bigskip
{\bf 4.2. Proposition.} {\it Assume that $M$ has Loewy length
at most $2$ and all the modules $\Omega^nM$ with $n\ge 1$ are bipartite.
Then $M$ is a Koszul module.

If $\soc{}_AA = J^2$, and $M$ is a Koszul module, then 
all the modules $\Omega^nM$ with $n\ge 1$ are bipartite.}
	\medskip
Proof. If $m\ge 0$ and $\Omega^{m+1}M$ is bipartite, then 
$\Omega^{m}M$ is aligned, see 3.4 (2).
Thus, the assumption implies that all the modules $\Omega^mM$ with $m\ge 0$ are aligned.
It follows from Proposition  4.1 that $M$ is Koszul. 

Assume now that $\soc{}_AA = J^2$. Then [RZ] 13.2 asserts that any aligned module
is bipartite. Thus, if $M$ is a Koszul module, then 
all the modules $\Omega^nM$ with $n\ge 1$ are aligned, thus bipartite.
$\s$
	\medskip
{\bf Remarks. (1)} If $A$ is a short local algebra with $a = 0$, 
then $\Omega S = S^e$ shows that 
$S$ is a Koszul module, thus $A$ is a Koszul algebra. 
	\smallskip
{\bf (2)} If $A$ is a Koszul algebra, then usually not all modules are Koszul. 
A typical example is  the $k$-algebra $A$ 
with generators $x,y$ and relations $x^2, xy, y^2$. Let $I = Ax \simeq A/Ax$.
As we see, $I$ is $\Omega$-periodic with period 1. It follows that $I$ is 
a Koszul module. 
But it is easy to see that the remaining indecomposable modules of length 2 are not
Koszul.

Also, all self-injective short local algebras $A$ with $e\ge 2$ are Koszul algebras,
but Proposition A.8 in the Appendix of [RZ] asserts that there are countably many indecomposable modules which are not Koszul.
	\medskip
{\bf (3)} Here is an example of 
{\it a Koszul module $M$ such that none of the modules $\Omega^nM$ with $n\ge 1$ 
is bipartite.} Let $A$ be  generated by $x,y$ with relations
$x^2, xy, y^2$ (this algebras has been considered already in [RZ], 9.3).
Let $I = Ax$. Then $\Omega S = {}_AJ = I\oplus S$, and 
$\Omega I = I.$ Thus, by induction, we see that $\Omega^n S = I^n\oplus S$
for all $n\ge 0.$ It follows that $M = S$ is a Koszul module. On the other hand,
$S$ is a proper direct summand of $\Omega^n M = \Omega^n S$, for any $n\ge 1.$
	\medskip

	\bigskip
{\bf 4.3.} We say that a short exact sequence $0 \to M' \to M \to M'' \to 0$
is {\it $t$-exact}, provided $t(M) = t(M')+t(M'')$. A
submodule $M'$ of $M$ will be called a
{\it $t$-submodule} provided the canonical exact sequence
$0 \to M' \to M \to M/M' \to 0$ is $t$-exact, thus provided $t(M) = t(M')+t(M/M').$
Of course, if $M'$ is a submodule of $M$, then $t(M) = t(M')+t(M/M')$
if and only if
$P(M)$ is isomorphic to $P(M')\oplus P(M/M').$

Similarly, 
a filtration $0 = M_0 \subseteq M_1 \subseteq \cdots \subseteq M_m = M$
will be called a {\it $t$-filtration} provided $t(M) = \sum_{j=1}^m t(M_j/M_{j-1}),$ 
or, equivalently, provided $P(M)$ is isomorphic to $\bigoplus_j P(M_j/M_{j-1})$.
Note that if $M$ has Loewy length at most 2 and 
$0 = M_0 \subseteq M_1 \subseteq \cdots \subseteq M_m = M$ is a $t$-filtration
then $\bdim M = \sum_j\bdim M_j/M_{j-1}.$
	\medskip
{\bf Lemma.} {\it Let $A$ be a short local algebra.
Let $0 \to M' \to M \to M' \to 0$ be $t$-exact. If $M', M''$ are
aligned, then also $M$ is aligned and there is a $t$-exact sequence
$0 \to \Omega M' \to \Omega M \to \Omega M'' \to 0.$}
	\medskip
Proof. We can assume that the map $M'\to M$ is the inclusion map of a submodule.
Since $P(M) \simeq P(M')\oplus P(M'')$, the horseshoe
lemma yields an exact sequence $0 \to \Omega M' \to \Omega M \to \Omega M'' \to 0.$
Multiplying with $J^2$, there is the exact sequence
$0 \to J^2\Omega M' \to J^2\Omega M \to J^2\Omega M'' \to 0$.
Now we have the inclusion maps $u'\:J\Omega M' \subseteq J^2\Omega M'$, 
$u\: J\Omega M \subseteq J^2\Omega M$, and $u''\:J\Omega M'' \subseteq J^2\Omega M''$.
If $M', M''$ are aligned, the maps $u',u''$ are bijective, thus also $u$ has
to be bijective. This shows that $M$ is aligned and therefore $\bdim \Omega M =
\omega(\bdim M).$ But this implies that 
$$
\align
 \bdim \Omega M &= \omega(\bdim M) = \omega(\bdim M'+\omega(\bdim M'') \cr
  &=
  \omega(\bdim M')+\omega(\bdim M'') 
  = \bdim \Omega M'+\bdim \Omega M''.
\endalign
$$
In particular, we have $t(\Omega M) = t(\Omega M')+t(\Omega M'').$ This shows that 
$\Omega M'$ can be identified with a $t$-submodule of $\Omega M$ with factor module $\Omega M''.$
$\s$
	\medskip
{\bf 4.4. Corollary.} {\it Let $A$ be a short local algebra. Let
$M$ be a module of Loewy length at most $2$ and
$0 = M_0 \subseteq M_1 \subseteq \cdots \subseteq M_m = M$ a $t$-filtration. Let $n\ge 0.$

{\rm(a)} If all the modules $\Omega^i(M_j/M_{j-1})$ with $0\le i \le n$ and $1 \le j \le m$
are aligned, then $\Omega^n M$ is aligned and $\Omega^{n+1}M$ has a $t$-filtration
with factors $\Omega^{n+1}(M_j/M_{j-1})$, where $1\le j \le m.$ 

{\rm(b)} If all the modules $M_j/M_{j-1}$ with $1\le j \le m$ are Koszul modules,
then also $M$ is a Koszul module.}
	\medskip 
Proof. Let $\omega = \omega^e_a$.
It is sufficient to show the assertions (a) and (b) for $m=2$, the general case follows easily by
induction on $m.$ Thus, let $M'$ be a submodule of $M$ and let $M'' = M/M'.$
We assume that $M'$ is a $t$-submodule of $M$, thus $t(M) = t(M')+t(M'')$ and  
$P(M) \simeq P(M')\oplus P(M'').$ 

	\smallskip
(a) We show:
{\it If $\Omega^iM', \Omega^iM''$ are aligned for $0 \le i \le n$, then $\Omega^n M$
is aligned and $\Omega^{n+1}M'$ can be identified with a $t$-submodule of $\Omega^{n+1}M$
with factor module $\Omega^{n+1}M''.$} Proof by induction on $n$. The case $n=0$ has been shown in 4.3.
Thus assume the assertion is true for some $n\ge 0$, and assume now that the modules
$\Omega^iM', \Omega^iM''$ are aligned for $0 \le i \le n+1$. By induction, we know that 
we may consider $\Omega^{n+1}M'$ as a $t$-submodule of $\Omega^{n+1}M$
with factor module $\Omega^{n+1}M''$. Since $\Omega^{n+1}M'$ and $\Omega^{n+1}M''$ are aligned, 
we apply (1) in order to conclude that $\Omega^{n+1}M$ is aligned and that 
$\Omega^{n+2} M'$ can be identified
with a $t$-submodule of $\Omega^{n+2} M$ with factor module $\Omega^{n+2} M''.$ This 
completes the proof of (a). 
	\smallskip
(b) {\it If $M',M''$ are Koszul modules, then $M$ is Koszul.} Proof. We use the equivalence of
(i) and (iv) in 4.1: If $M',M''$ are Koszul, then
all the modules $\Omega^n M', \Omega^n M''$
are aligned. According to (a), all the modules $\Omega^nM$ are aligned. Thus $M$ is Koszul.
$\s$
	\bigskip
{\bf 4.5. Corollary.} {\it Let $A$ be a short local algebra and $U$ an ideal of $A$. Let $n\ge 0.$

{\rm (a)} Assume that for any local module $N$ annihilated by $U$, the modules $\Omega^iN$ with
$0 \le i \le n$ are aligned. Then for any module $M$ annihilated by $U$, 
the module $\Omega^nM$ is aligned.

{\rm (b)} Assume that any local module $N$ annihilated by $U$ is a Koszul module,
then any module $M$ annihilated by $U$ is a Koszul module.}
	\medskip
Proof. Any module $M$ annihilated by $U$ has a $t$-filtration whose factors are local modules
(of course annihilated by $U$). Namely, any composition series of $\top M$ lifts to
a $t$-filtration of $M$. Thus, we can apply 4.4.
$\s$
	\medskip
{\bf 4.6. Proposition.} {\it Let $A$ be a short local algebra. If there exists 
a non-projective Koszul module, then $A$ is a left Koszul algebra.}
	\medskip
Proof. If $N$ is a non-projective Koszul module, then $\Omega N$ is a non-zero
module of Loewy length
at most $2$ and is a Koszul module. Thus we can assume that there is given a module
$M\neq 0$ of Loewy length at most 2 which is a Koszul module. 
According to 3.2, we have the exact sequence
$$
 \eta_M = (\ 0 @>>> \Omega M @>u>> JP(M) @>p>> JM @>>> 0\ ).
$$
Since $M$ is a Koszul module, $M$ is aligned, thus 3.3 (vii) asserts that $\eta_M$
is $t$-exact.

Consider the sequences
$$
 \eta_M(n) =(\ 0 @>>> \Omega^{n+1} M 
  @>\Omega^nu>> \Omega^n(JP(M)) @>\Omega^n p>> \Omega^n(JM) @>>> 0\ ).
$$
with $n\ge 0.$ By induction on $n$ we show that $\eta_M(n)$ is $t$-exact, 
and that all modules $\Omega^{n}S$ is aligned. 

Proof of the induction. First, let $n = 0$. We know that $\eta_M(0) = \eta_M$ is
$t$-exact. Also, the module 
$\Omega^0 S = S$ is always aligned. Now assume that for some $n\ge 0$
the sequences $\eta_M(n)$ is $t$-exact and the modules $\Omega^n S$ 
is aligned. Since $M$ is 
a Koszul module, the module $\Omega^{n}\Omega M$ 
is aligned. Since $JM$ is semisimple, 
the modules $\Omega^{n}(JM)$ is aligned.
Thus, we can apply Lemma 4.3  in order to conclude that
$\Omega^n(JP(M))$ is aligned and that $\eta_M(n+1)$ is $t$-exact.
Let $m = |\top M|$. Then $JP(M) = \Omega(S^m)$. Since
$\Omega^{n} JP(M) = \Omega^{n+1} S^m$ is aligned and $m\ge 1$, we see that
$\Omega^{n+1} S$ is aligned. This completes the induction step.

Altogether, we see that $\Omega^nS$ is aligned for all $n\ge 0$, thus
$S$ is a Koszul module.
$\s$
	\bigskip
{\bf 4.7.} Finally, let us draw the attention again to the simple module $S$.
By definition, $A$ is a left Koszul algebra iff $S$ is Koszul. 
What does it mean that $S$ is a Koszul module?

If $e, a $ are real number, one may define recursively 
the sequence $b_n = b(e,a)_n$ with $n\ge -1$ as follows:
$b_{-1} = 0,$ $b_0 = 1$ and 
$$
  b_{n+1} = eb_{n}-ab_{n-1}, \tag{$*$}
$$
for $n\ge 0.$ By induction, one sees that $(b_n,ab_{n-1}) = (\omega^e_a)^n(1,0).$ 
	\medskip
{\bf Proposition.} {\it Let $A$ be a short local algebra of Hilbert type $(e,a)$.
The module $S$ is Koszul iff 
$\bdim \Omega^n S = (b(e,a)_n,a\cdot b(e.a)_{n-1})$  for all $n\ge 0$.}
	\medskip 
Proof. Write $b_n = b(e,a)_n$ for all $n\ge -1.$ 
According to 4.1, $S$ is a Koszul module iff $\bdim \Omega^n S 
= (\omega^e_a)^n \bdim S$ for all $n\ge 0.$ Of course, $\bdim S = (1,0) = (b_0,b_{-1}),$ 
and therefore
$(\omega^e_a)^n \bdim S  = (\omega^e_a)^n(1,0) = (b_n,ab_{n-1}).$
$\s$
	\medskip
{\bf Remark.}
Avramov-Iyengar-\c Sega have shown: if  $a < \frac14e^2$, then 
for all $n\ge 0$
$$
 b(e,a)_n =  \frac 1{2^n} \sum_{j=0}^{\lfloor \frac {n}2\rfloor} 
 \binom{n+1}{2j+1} (e^2-4a)^{j}e^{n-2j},
$$
see Appendix B of [RZ]. 

	\bigskip\bigskip
\vfill\eject
{\bf 5. Again: The $\Omega$-growth of a module.}
	\medskip
Let $A$ be a short local algebra and $M$ a module of Loewy length at most $2$.
What are the possible values for $\gamma(M)$? First, we assume that $M$ is not
a Koszul module. The following observation was mentioned already in the introduction,
see 1.3. 
	\medskip
{\bf 5.1. Proposition.} {\it Let $A$ be a short local algebra. Let $M$ be a module
of Loewy length at most $2$.
If $M$ is not Koszul, then $\gamma(M) = \gamma_A$.}
	\medskip
Proof. Proposition 4.2 asserts that there is $n\ge 1$ such that $\Omega^n M$
is not bipartite. Since $\Omega^n M$ has Loewy length at most 2, we see that
$S$ is a direct summand of $\Omega^n M$. According to Theorem 1,  
$\gamma(M) = \gamma_A.$
$\s$
	\bigskip
It remains to consider the Koszul modules.
We will need two elementary considerations from real linear algebra.
Given vectors $\bold x = (x_1,x_2)$ and $\bold y = (y_1,y_2) \in \Bbb R^2,$ 
we write $\bold x \le \bold y$ provided $x_1\le y_1$ and $x_2 \le y_2$
and we write 
$\bold x < \bold y$ provided $\bold x \le \bold y$ and $\bold x \neq \bold y$.
Let $|\bold x| = |x_1|+|x_2|.$

If $\omega\:\Bbb R^2 \to \Bbb R^2$ is a linear transformation, let $\gamma_\omega(\bold x)
= \limsup_n \root n \of {|\omega^n(\bold x)|}.$
	\medskip
{\bf 5.2. Lemma.} {\it Let $\omega\:\Bbb R^2 \to \Bbb R^2$ be a linear transformation. 
If $\bold x$ is an eigenvector of $\omega$ with eigenvalue $\lambda$, then
$\gamma_\omega(\bold x) = |\lambda|.$
If $\bold x$ is non-zero and not an eigenvector of $\omega$, 
then $\gamma_\omega(\bold x) = \rho(\omega)$.}
	\medskip

Proof. Of course, if $\omega(\bold x) = \lambda\bold x$, then 
$\gamma_\omega(\bold x) = |\lambda|.$ Thus, let us assume that $\bold x = (x_1,x_2)$ 
is non-zero and not an eigenvector of $\omega.$ 

First, let $\omega$ be semisimple with eigenvalues $\lambda, \lambda',$ where 
$|\lambda'|\le |\lambda|.$ Thus we can assume that $\omega = \left[\smallmatrix \lambda & 0\cr
                         0 & \lambda' \endsmallmatrix\right]$.  Then 
$\omega^n(x_1,x_2) = \lambda^n(x_1+(\lambda'/\lambda)^n x_2)$.
Since $0 \le |(\lambda'/\lambda)^n x_2)| \le |x_2|$ and $x_1\neq 0,$ we see
that 
$\gamma_\omega(x_1,x_2) = |\lambda|\cdot
\limsup_n\root n\of {|x_1|+ |(\lambda'/\lambda)^n x_2|} =
|\lambda| = \rho(\omega).$

Second, let $\omega$ be not semisimple. Let $\lambda$ be its eigenvalue. If $\lambda = 0,$
then $\gamma_\omega(\bold x) = 0 = \rho(\omega).$
Otherwise, we can assume that $\omega = \lambda \left[\smallmatrix 1 & 1\cr
                         0 & 1 \endsmallmatrix\right]$,
thus $\omega^n(x_1,x_2) = \lambda^n(x_1+nx_2,x_2)$ and 
$\gamma_\omega(x_1,x_2) = 
|\lambda|\cdot\limsup_n\root n\of {|x_1+nx_2|+|x_2|} = |\lambda| = \rho(\omega).$
$\s$
	\medskip

{\bf 5.3. Lemma.} {\it Let $\omega\:\Bbb R^2 \to \Bbb R^2$ be a linear transformation.
 Let $\bold x$ be an element of $\Bbb R^2$ such that $\omega^n \bold x > 0$
for all $n\ge 0$. Then $\omega$ has real eigenvalues. If $\bold x$ is an eigenvector with
eigenvalue $\lambda$, then $\lambda$ is positive. If $\bold x$ is not an eigenvector, then
$\rho(\omega)$ is an eigenvalue of $\omega$ (and, of course, positive).}
	\medskip
Note that this lemma is a version of the Perron-Frobenius theorem in dimension 2, but
in contrast to the classical Perron-Frobenius theorem, we cannot assert hat
$\gamma_\omega(\bold x)$ is a simple eigenvalue of $\omega$, as the example of $\omega = 
\left[\smallmatrix 2 & -1\cr
                   1 & 0 \endsmallmatrix\right]$
and $\bold x = (1,0)$ shows: We have 
$\omega^n \bold x = (n+1,n)> 0$ for all $n\ge 0$, and $\bold x$ is not an eigenvector
of $\omega$; on the other hand, $1$ is an eigenvalue of $\omega$ with multiplicity 2.
	\medskip
Proof of Lemma. 
Let $\rho = \rho(\omega)$ be the  spectral radius of $\omega$.
Let $\bold x$ be a vector in $\Bbb R^2$ with $\omega^n \bold x > 0$
for all $n\ge 0$. The existence of $\bold x$ shows that $\omega$ cannot be nilpotent, thus 
$\rho(\omega) > 0.$ Of course, if $\bold x$ is an eigenvector with eigenvalue $\lambda$, then
$\lambda > 0.$ Thus, let us assume that $\bold x$ is not an eigenvector. 

Given a set $\Cal X$ of vectors in $\Bbb R^2$, let $C(\Cal X)$ 
be the cone in $\Bbb R^2$ of all vectors which are
linear combinations of the elements in $\Cal X$  using positive coefficients.
Let $C = C(\{\omega^n(\bold x)\mid n\ge 0\})$.
Then all non-zero vectors $\bold y\in C$ satisfy $\bold y > 0$ and $\omega(C) \subseteq C$. 
If $C$ is a ray, then any non-zero element in $C$ is an eigenvector of $\omega$, thus
$\bold x$ is an eigenvector, a contradiction.

Thus, $C$ is not a ray, and 
there is a basis $\bold y,\ \bold y'$ of $\Bbb R^2$
such that the topological closure $\overline C$
of $C$ is the cone $\overline C = C(\{\bold y,\bold y'\})$.
We have $\omega(\overline C) \subseteq \overline C$, in particular 
$\omega(\bold y), \omega(\bold y')\in 
\overline C.$ 

If $\omega(\bold y) \in \Bbb R_+\bold y$, say $\omega(\bold y) = \lambda\bold y$
with $\lambda \in \Bbb R_+,$ let $\omega(\bold y') = c\bold y + d\bold y'$, thus
$c\ge 0$ and $d>0.$ Now $\omega$ is similar to the matrix
$\left[\smallmatrix \lambda & c\cr
                                0 & d \endsmallmatrix\right]$, thus its eigenvalues are
$\lambda$ and $d$, and both are positive. Therefore $\rho(\omega) = \max\{\lambda,d\}$
is an eigenvalue.

Next, assume that $\omega(\bold y) = \lambda\bold y'$. Let $\bold y' = c\bold y+d\bold y'$.
Then $\omega$ is similar to a matrix of the 
form $\left[\smallmatrix 0 & c\cr
                         1 & d \endsmallmatrix\right]$,
its eigenvalues are $\frac12d\pm \frac12\sqrt{d^2+4c}$, thus the spectral radius is
the eigenvalue $\frac12d + \frac12\sqrt{d^2+4c}$.

Finally. it remains to consider the case that both $\omega(\bold y), \omega(\bold y')$ 
belong to the interior of $\overline C$, then
$\omega$ is similar to a matrix with positive coefficients, thus the usual Perron-Frobenius
theorem asserts that the spectral radius of $\rho(\omega)$ is an  eigenvalue
of $\omega$.
$\s$
	\medskip
{\bf 5.4. Lemma.} {\it The transformation
$\omega^e_a$ has real eigenvalues iff $a \le \frac14e^2$. In this case, both
eigenvalues are non-negative and $\rho(\omega^e_a) = \frac12(e + \sqrt{e^2-4a}).$

If $\lambda\neq 0$ is an eigenvalue of $\omega^e_a$, then $(\lambda,a)$ 
is eigenvector of $\omega^e_a$ with eigenvalue $\lambda$.}
	\medskip
Proof: The eigenvalues of $\omega^e_a$ are $\frac12(e\pm \sqrt{e^2-4a})$, thus they are real iff
$e^2 \ge 4a.$ Also, since $e^2 - 4a \le e^2,$ it follows from $e^2\ge 4a$ that 
$\sqrt{e^2-4a} \le e.$ 

Let $\lambda$ be an eigenvalue of $\omega^e_a$. The characteristic polynomial
of $\omega^e_a$ is $T^2-eT+a$, thus $\lambda^2 = e\lambda - a$, therefore
$\omega^e_a(\lambda,a) = (e\lambda-a,a\lambda) = \lambda(\lambda,a).$ 
$\s$

	\bigskip 
{\bf 5.5. Proposition.} 
{\it Assume that $A$ is a short local algebra of Hilbert type $(e,a)$
with $e\ge 2$. If there exists a non-projective Koszul module $M$, then $a \le \frac14e^2$
and $\gamma(M)$ is a positive eigenvalue of $\omega^e_a.$}
	\medskip
Proof. Let $M$ be a non-projective Koszul module. Replacing, if necessary, $M$ by $\Omega(M)$,
we can assume that $M$ has Loewy length at most 2. Let 
$\bold x = \bdim M.$ Since $M$ is Koszul, we have $(\omega^e_a)^n \bold x = \bdim \Omega^n M > 0$
for all $n\ge 0.$ Thus, 5.3 
assert that $\gamma(M)$ is a positive real eigenvalue. 
The existence of a real
eigenvalue shows that $a\le \frac14e^2$, see 5.4.
$\s$
	\medskip
{\bf Corollary.} {\it Let $A$ be a short local algebra of Hilbert type $(e,a)$ with
$a > \frac14e^2$. If $M$ is a non-projective module, then $M$ is not Koszul,
thus $\gamma(M) = \gamma_A$.}
	\medskip
Proof. Let $M$ be non-projective. Since $a > \frac14e^2$, 
Proposition 5.5 asserts that $M$ cannot be Koszul.
According to 5.1 we have $\gamma(M) = \gamma_A.$ 
	\bigskip
{\bf 5.6. Proof of Theorems 2 and 3.} Let $A$ be a short local algebra of Hilbert type $(e,a)$.

First, let $a = 0$. Then $\Omega^S = S^e$ shows that $S$ is a Koszul module 
($S$ is always aligned) and $\gamma_A = e$. 
For any non-zero module $M$, the module $\Omega M$ is semisimple and not zero,
thus $\gamma(M) = \gamma_A$.

Now let $a \neq 0.$ If $e = 1$, then $\Omega J = S$, thus 
$\Omega^2S = S$ 
shows that $S$ is not Koszul, thus there are no non-projective Koszul modules. 

Thus, let $e \ge 2$ and $a\neq 0.$  
Let $M\neq 0$ be a Koszul module of Loewy length at most 2.
According to 4.6, $A$ is left Koszul. According to 5.5 we know that $a\le \frac14 e^2$
and that $\gamma(M)$ and $\gamma_A$ are positive eigenvalues of $\omega^e_a.$ 
Since $a \neq 0$, we see that $(1,0)$ is not an eigenvector of $\omega^e_a$, thus
5.2 asserts that $\gamma_A = \gamma_\omega(1,0) = \rho(\omega) = 
\frac12(e+\sqrt{e^2-4a}\,)$. 
Assume that $\gamma(M) \neq \gamma_A$, then $\gamma(M) \neq \rho(\omega)$,
thus 5.2 asserts that $\bdim M$ is an eigenvector of $\omega$
and $\gamma(M)$ is the corresponding eigenvalue, thus equal to
$\frac12(e-\sqrt{e^2-4a}\,)$. This completes the proof of Theorem 2.
	\smallskip
Now assume that there is a non-zero module $M$ of Loewy length at most 2 with
$\gamma(M) 
< \gamma_A$. As we have seen, $\bdim M$ is an eigenvector of $\omega$ and the
corresponding eigenvalue is $\gamma(M)$. But this means that $\gamma(M)\bdim M$
is a vector with integral coefficients, thus $\gamma(M)$ has to be rational
and therefore $e^2-4a$ has to be the square of an integer. 
Since $\bdim M$ is an eigenvector of $\omega$ with eigenvalue $\gamma(M)$,
and all eigenvectors have multiplicity 1, 
$\bdim M$ is a multiple of $(\gamma(M),a) = (\gamma(M),\gamma(M)\gamma_A$, and thus
a multiple of $(1,\gamma_A).$
	\smallskip
Let us assume that $e =  \gamma(M)+\gamma_A$.
Since $0 < \gamma(M) < \gamma_A$, we have $2\gamma(M) < \gamma(M)+\gamma_A = e,$
thus $\gamma(M) < \frac12e.$ Since $\gamma(M) < \frac12e$, we have $\frac12 e = e-\frac12e < 
e-\gamma(M) = \gamma_A.$ Since $0 < \gamma(M)$, we have $\gamma_A < 
\gamma(M)+\gamma_A = e.$ This shows that $0 < \gamma(M) < \frac12e < \gamma_A < e$.
Since $S$ is a Koszul module, Theorem 2 asserts that $\gamma_A = \frac12(e+\sqrt{e^2-4a}\,)$
and $\gamma(M) = \frac12(e-\sqrt{e^2-4a}\,).$
It follows that $\gamma_A-\gamma(M) = \sqrt{e^2-4a},$ thus
$(\gamma_A -\gamma(M))^2 = e^2-4a$, so that $e^2-4a$ is the square of a positive integer.
$\s$
	\bigskip
{\bf 6. Left Conca ideals.}
	\medskip
%%%%%%%%%%%%%%%%%%%%%%%%%%%%%%%%%%%%%%%%%%%%%%%%%%%%%%%%%%%%%
{\bf 6.1.} 
Let $A$ be a local algebra and $U$ an ideal of $A$. We say that $U$ is a
{\it left Conca} ideal provided $U^2 = 0$ and $J^2 \subseteq JU$.
If $A$ has a left Conca ideal $U$, then $A$ 
is short (namely, $J^3 \subseteq J^2U \subseteq JU^2 = 0$).
	\medskip 
{\bf Remark.} The name corresponds to the considerations in 
[AIS]. Following [AIS] (but dealing also with non-commutative local algebras),
an element $x$ may be called a 
{\it left Conca generator} of $J$ provided $x^2 = 0 \neq x$ and $J^2 = Jx$.
If $x$ is a left Conca generator of $J$, then clearly $Ax$ is a left Conca ideal
(note that $Ax$ is a twosided ideal, since $xA = kx + xJ \subseteq kx + J^2 \subseteq
Ax$). 
Obviously, the existence of a left Conca generator for $J$ implies that $a \le e-1.$
As we will see in 7.1, for any pair $(e,a)$ with $a\le \frac14 e^2$, there are
short local algebras of Hilbert type $(e,a)$ with a left Conca ideal. 
	\medskip
{\bf 6.2. Proof of Theorem 4.}  Let $U$ be a left Conca ideal in $A$. 
Let $N$ be a local module annihilated by $U$, 
thus, $N \simeq {}_A/V$ for some proper left ideal $V$ of $A$ and $U \subseteq V$, since
$UN = 0.$ Since $J^2 \subseteq JU \subseteq U \subseteq V \subseteq J$, the factor module
$V/U$ is a subquotient of ${}_AJ/J^2$, thus semisimple. in addition, 
$J^2 \subseteq JU \subseteq JV \subseteq J^2$ shows that the embedding $u\:U \to V$ yields
the equality $JU = JV$. Thus, the Snake Lemma applied to
$$
\CD 
 0 @>>> JU @>1>> JV @>>> 0 @>>> 0 \cr
 @.     @VVV  @VVV      @| \cr
 0 @>>> U @>u>>  V @>>> V/U @>>> 0
\endCD
$$
shows that $U$ is a $t$-submodule of $V$. 

We show by induction on $n\ge 0$: {\it If $N$ is a local module annihilated by $U$, then
$\Omega^i N$ is aligned, for $0\le i\le n.$}

First, let $n = 0.$ We have $N = {}_AA/V$ for some left ideal $V$, thus
$\Omega N = V$, and as we have mentioned already,
$J^2 \subseteq JV = J\Omega N$, thus condition (v) of 3.3 asserts that $N$ is aligned. 

Now, assume that we know for some $n\ge 0,$ that for all local modules $N'$ annihilated by $U$
the modules $\Omega^i N'$ with $0\le i \le n$ are aligned. According to Corollary (a) in
4.5, this implies that for all modules $M$ annihilated by $U$, the modules
$\Omega^i M$ with $0\le i \le n$ are aligned. Let $N$ be a local module annihilated by $U$,
say $N = {}_AA/V$ for some left ideal $V$. By induction assumption, we know that the modules
$\Omega^iN$ are aligned for $0\le i \le n.$ It remains to be seen that $\Omega^{n+1}N$ is
aligned. As we have mentioned, $U$ is a $t$-submodule of
$V$. Now $U$ is annihilated by $U$. Also, $V/U$ is annihilated by $U$ (since it is semisimple).
Thus, all the modules
$\Omega^iU$ and $\Omega^i(V/U)$ are aligned, for $0\le i \le n$. We apply Lemma 4.4 (a) 
in order to conclude that $\Omega^n V$ is aligned. Thus, $\Omega^{n+1}N = \Omega^nV$ is aligned. 
$\s$
	\medskip 
{\bf Remark.} This improves Theorem 3.2 of [AIS]. In addition, we should stress
that the proof yields the following stronger assertion: 
{\it If $A$ has a left Conca ideal $U$, then any module 
with a $t$-filtration with factors annihilated by $U$ is a Koszul module.}

One should be aware that given any ideal $U$, there may be modules which are not
annihilated by $U$, but which have a $t$-filtration
with factors annihilated by $U$. 
For example, if $A$ is of Hilbert type $(2,0)$ (thus, $A$ is the local 3-dimensional
algebra with radical square zero)
and $U$ is one-dimensional, then there are
just two indecomposable modules annihilated by $U$, namely $S$ and $I = {}_AA/U$, 
but infinitely many indecomposable modules 
which have a $t$-filtration with factors of the form $I$ and $S$, namely the modules in the  Auslander-Reiten component which contains $I$ 
as well as the preinjective modules.
	\medskip
{\bf 6.3. Remark.}
{\it A short local Koszul algebra $A$ 
may not have any left Conca ideal. Also, $A$ may not have a 
left Conca ideal, whereas its opposite algebra has a left Conga ideal.} 

Here is an example:
Let $A$ be generated by $x,y,z$ with relations 
$$ 
   x^2,\ yx,\ zx,\ zy,\ y^2-xz,\ yz,\ z^2,
$$
so that $J^2$ has the basis $xy,\ y^2=xz.$ One easily checks that $A$ has no left Conga ideal (namely, any ideal $U$ with $U^2 = 0$
is contained in $Ax+Az$, thus $JU \subseteq kxz$).
But $xA$ (with basis $x,\ xy,\ xz$ is a right Conga ideal. 

Since the opposite algebra of $A$ has a left Conga ideal, $A$ is a right
Koszul algebra. In order to see that $A$ is also 
left Koszul, write ${}_AJ = S\oplus W$ where $W = Ay+Az$. Then
$\Omega W = W^2$. Therefore $\Omega^nS = S\oplus W^n$ has dimension vector 
$(2n+1,2n) = (\omega^3_2)^n (1,0)$ and this shows that $S$ is a Koszul module. 
$\s$
	\bigskip
%%%%%%%%%%%%%%%%%%%%%%%%%%%%%%%%%%%%%%%%%%%%%%%%%%%%%%%%%%%%%
{\bf 7. Construction of Koszul algebras.}
	\medskip
{\bf 7.1. Proposition.} {\it 
If $0\le  a\le \frac14e^2$, there are short local algebras of 
Hilbert type $(e,a)$ (even commutative ones) with a left Conca ideal.}

	\medskip
Proof. Assume that 
$0\le a\le \frac14e^2.$ We are going to construct a commutative short local algebra $A$
of Hilbert type $(e,a)$ which is Koszul. 

Let $c = \lfloor\frac12e\rfloor$ and $d = e-c$. Since $0\le a\le \frac14e^2$, we have
$a\le cd$ (namely, for $e$ even, $c = d = \frac12e$ and $a\le c^2 = cd,$  whereas for
$e$ odd, we have $d = c+1$ and $a \le \frac14(2c+1)^2$ implies that $a\le c^2+c = cd$).
Thus we can write $a = \sum_{j=1}^d a(j)$ with $0\le a(j)\le c.$ 

Let $A = \Lambda(c;a(1),\dots,a(d))$ be the commutative algebra generated by the 
elements $x_i,\ y_j$ 
with $1\le i \le c$ and $1\le j \le d$ and the relations
$x_ix_{i'}, y_jy_{j'}$ for all $i,i'\in\{1,\dots,c\}$ and
$j,j'\in\{1,\dots,d\},$ as well as $x_iy_j$ for all pairs $(i,j)$
with $1\le j \le d$ and $a(j) < i \le c.$ 
It follows that $J^2$ has the basis $x_iy_j$ with $1\le j \le d$ and
$1\le i \le a(j).$ 

If $U = \sum_{j=1}^d Ay_j,$ then $U^2 = 0$ and $J^2 \subseteq JU,$ thus $U$ is a left Conca ideal. 
$\s$
	\bigskip
{\bf 7.2. Proof of Theorem 5.} If $A$ is a short local left Koszul algebra of 
Hilbert type $(e,a)$, then Theorem 2 asserts that $0 \le a \le \frac14 e^2.$ 
Conversely, 7.1 shows that for $0\le  a\le \frac14e^2$, there are commutative
short local algebras $A$ of Hilbert type $(e,a)$ with a left Conca ideal.
According to Theorem 4, these algebras $A$ are left Koszul algebras.
$\s$

	\bigskip
%%%%%%%%%%%%%%%%%%%%%%%%%%%%%%%%%%%%%%%%%%%%%%%%%%%%%%%%
{\bf 7.3.}
{\it For any pair $c,d$ of natural numbers, there exists a commutative
short local algebra $\Lambda(c,d)$
of Hilbert type $(e,a)$, where $e = c+d$ and $a = cd$,
with a module $M$ with dimension vector $(1,c)$ such that
$\Omega M \simeq M^d$ (thus $\gamma(M) = d$), and such that $\gamma_A = \max\{c,d\} 
= \rho(\omega^e_a).$} 
	\medskip 

Proof. Let $A = \Lambda(c,d)$ be the commutative algebra
generated by $x_1,\dots,x_c,\ y_1,\dots,y_d,$ 
and with relations $x_ix_{i'}, y_jy_{j'}$ for all $i,i'\in\{1,\dots,c\}$ and
$j,j'\in\{1,\dots,d\}.$ Then $J^2$ has the basis $x_iy_j$ with $1\le i \le c,\ 1\le j \le d.$
Let $M = A y_1$; this is a local module of Loewy length $2$ with socle $x_1y_1,\dots,x_cy_1$,
thus with dimension vector $(1,c)$. All the module $A y_j$ with $1\le j \le d$
are isomorphic to $M$ and
$J$ is isomorphic to $S^c\oplus M^d.$ Since $A/\bigoplus_{j=1}^d A y_j$ is isomorphic to $M$,
we see that $\Omega M \simeq M^d.$ It follows that $\gamma(M) = d.$ 

Here is a similar, but non-commutative example: {\it a non-commutative short local algebra 
$\Lambda'(c,d)$
of the same Hilbert type $(c+d,cd)$ with a module $M$ with dimension vector $(1,c)$ such that
$\Omega M \simeq M^d,$ so that $\gamma(M) = d,$ whereas $\gamma_A = \max\{c,d\} = \rho(\omega^e_a).$}
Let $\Lambda'(c,d)$ be 
generated by $x_1,\dots,x_c,\ y_1,\dots,y_d,$ 
and with relations $x_ix_{i'}, y_jy_{j'}, y_jx_i$ for all $i,i'\in\{1,\dots,c\}$ and
$j,j'\in\{1,\dots,d\}.$ Again, $J^2$ has the basis $x_iy_j$ with $1\le i \le c,\ 1\le j \le d$.
Note that the elements $x_1,\dots,x_d$ do not belong to $J^2$, but to 
$\soc{}_{\Lambda'(c,d)}J$. 
	\medskip
{\bf 7.4.}
In particular, let us focus the attention to the case $d=1.$ 
The algebras $\Lambda(c,1)$ and $\Lambda'(c,1)$ have a non-zero $\Omega$-periodic 
module $M$.

On the other hand, let us stress that 
the algebra $A = \Lambda'(c,1)$ is a short local algebra with $J^2 \subset \soc{}_AA$
as well as $J^2 \subset \soc A_A$. Note that 
Lescot [L] Prop. 3.9 (2) has pointed out that for a {\bf commutative} short local algebra with
$J^2 \subset \soc A$ and a non-projective module $M$, the sequence $t_n(M)$ is always strictly
increasing. 
	\bigskip\bigskip 
%%%%%%%%%%%%%%%%%%%%%%%%%%%%%%%%%%%%%%%%%%%%%%%%%%%%%%%%%%%%%
{\bf 8. A lower bound for $\gamma_A$.}
	\medskip
{\bf 8.1. Proof of Theorem 7.} We assume that $A$ is a short local algebra of Hilbert type
$(e,a)$ with $a\le \frac 14 e^2.$ Let $\omega = \omega^e_a$ and $\bold d(n) = \omega^n(1,0)$
for all $n\ge 0.$ We have $\omega(0,-1) = (1,0) = \bold d(0),$ and therefore
$\omega(w,-w) = w\bold d(1)+w\bold d(0)$ for any $w\in \Bbb Z.$  
	\smallskip
(1) Let us show that $\bold d(n) > 0$ and that 
$$
 \limsup_n \root n \of {|\bold d(n)|} = \tfrac12 (e+\sqrt{e^2-4ac}\,). 
$$
According to Theorem 5, there exists a short local algebra $A'$ of Hilbert type $(e,a)$
which is left Koszul. Let $S'$ be the simple $A'$-module. Since $S'$ is a Koszul module,
we have $\bdim \Omega_{A'}S' = \bold d(n),$ thus $\bold d(n) > 0.$
Theorem 2 asserts that
$\limsup_n \root n \of {|\bold d(n)|} = \limsup_n \root n \of {|\Omega_{A'}^nS'|} =
\frac12 (e+\sqrt{e^2-4ac}\,)$. 
	\smallskip
Let $\Cal N(n)$ be the set of linear combinations of $\bold d(i)$ with $0\le i \le n-1$ 
using non-negative coefficients. For $n\ge 1$, we apply the Main Lemma 3.1 to $\Omega^{n-1}S$
and obtain $\bdim \Omega^n S = \omega (\bdim \Omega^{n-1}S) + (w_n,-w_n)$ 
for some integer $w_n \ge 0.$ In addition, we define
$w_0 = 0.$ 
	\smallskip
(2) Using induction on $n\ge 0$, we show that 
$$
 \dim \Omega^n S - \bold d(n) - (w_n,-w_n) \ \in \ \Cal N(n).
$$

Proof. The assertion holds true for $n = 0$,
since $\bdim S = \bold d(n)$ and $w_0 = 0.$
Now assume that the assertion is true for some $n\ge 0$, thus we have
$$
 \dim \Omega^n S = \bold d(n) + \bold x \quad \text{with} \quad
  \bold x = (w_n,-w_n) + \sum_{i=0}^{n-1} v_i\bold d(i)
$$
with non-negative integers $v_i$, where $0\le i < n$.

We apply $\omega$ to $\bold x$
and get
$$
\align
  \omega(\bold x) &= \omega(w_n,-w_n) + \sum_{i=0}^{n-1} v_i\omega(\bold d(i)) \cr
                  &= w_n\bold d(1) + w_n\bold d(0) +\sum_{i=0}^{n-1} v_i\bold d(i+1), 
\endalign
$$
thus $\omega(\bold x)$ belongs to $\Cal N(n+1)$. On the other hand, we have 
$$
\align
 \bdim \Omega^{n+1} S &= \omega(\bdim\Omega^{n} S) + (w_{n+1},-w_{n+1}) \cr
    &= \omega(\bold d(n)+\bold x) + (w_{n+1},-w_{n+1}) \cr
    &= \bold d(n+1)+ (w_{n+1},-w_{n+1}) +\omega(\bold x). \cr
\endalign
$$
This shows that $\bdim \Omega^{n+1} S -  \bold d(n+1) - (w_{n+1},-w_{n+1}) = \omega(\bold x)$,
and we have seen already that $\omega(\bold x)$ 
belongs to $\Cal N(n+1).$ $\s$
	\smallskip
(3) We have $|\Omega^n S| \ge |\bold d(n)|$ for $n\ge 0.$ 
Namely, 
the formula $(2)$ implies that $\dim \Omega^n S  \ge \bold d(n) + (w_n,-w_n)$
(since $\bold d(i) \ge 0$ for all $0\le i < n$)
and therefore $ |\Omega^n S| \ge |\bold d(n)|,$
since $|(w_n,-w_n)| = 0.$ 
	\smallskip
(4) Altogether, (1) and (3) show that 
$$
 \gamma_A = 
 \limsup_n \root n \of {|\Omega^n S|} \ge 
 \limsup_n \root n \of {|\bold d(n)|} = \tfrac12 (e+\sqrt{e^2-4ac}\,),
$$
this completes the proof.
$\s$
	\bigskip
{\bf 8.2.} Let us show that $\gamma_A$ does not only depend on the Hilbert type. Of course,
as we have seen, if $A$ is a Koszul algebra, then $\gamma_A$ is determined by the
Hilbert type $(e,a)$, namely $\gamma_A = \rho(\omega^e_a)$. But we will show that there are
algebras $A,A'$ which are not Koszul with $\gamma_A \neq \gamma_{A'}$ (and both $\gamma_A,
\gamma_{A'}$ different from $\rho(\omega^e_a)$).
	\medskip
{\bf Example.} {\it 
Short local algebras $A$ of Hilbert type $(3,2)$ with $\gamma_A = 2,\ \psi,\ 3$, 
where $\psi = \frac12(3+\sqrt5)$ is the square of the golden ratio.} Note that $\rho(\omega^3_2) =
2$. 
	\medskip
First. If a short local algebra
has Hilbert type $(3,2)$ and is Koszul, then theorem 2 asserts that
$\gamma(S) = \rho(\omega^e_a) =
\frac12(3+\sqrt{9-4\cdot 2}\,) = 2$ (and we know from Theorem 5 that such algebras
do exist).
	\medskip
We define two algebras $A, A'$ of Hilbert type $(3,2)$ with generators $x,y,z$.
The relations for $A$ are $yx,\ zx,\ y^2,\ zy,\ xz,\ yz,\ z^2.$
The relations for $A'$ are $yx,\ zx,\ xy,\ zy,\ xz,\ yz,\ z^2.$
The radicals $J$ and $J'$, respectively, look as follows:
$$
{\beginpicture
    \setcoordinatesystem units <.6cm,1cm>
%%%%%%%%%%%%%%%%%%%%%%%%%%%%%%%%%%%%%%%%%%
\put{\beginpicture
\put{$x$\strut} at 0 1
\put{$y$\strut} at 2 1
\put{$z$\strut} at 4 1
\put{$x^2$\strut} at 1 0
\put{$zy$\strut} at 3 0
\arr{0.2 0.8}{0.8 0.2}
\put{$\ss x$\strut} at 0.15 0.5 
\put{$\ss z$\strut} at 2.15 0.5 
\arr{2.2 0.8}{2.8 0.2}
\put{$J$} at -1.5 1
\endpicture} at 0 0 
%%%%%%%%%%%%%%%%%%%%%%%%%%%%%%%%%%%%%%%%%%
\put{\beginpicture
\put{$x$\strut} at 0 1
\put{$y$\strut} at 2 1
\put{$z$\strut} at 4 1
\put{$x^2$\strut} at 1 0
\put{$y^2$\strut} at 3 0
\arr{0.2 0.8}{0.8 0.2}
\put{$\ss x$\strut} at 0.15 0.5 
\put{$\ss y$\strut} at 2.15 0.5 
\arr{2.2 0.8}{2.8 0.2}
\put{$J'$} at -1.5 1
\endpicture} at 9 0 
\endpicture}
$$
We will show that $\gamma_A = \psi$ and $\gamma_{A'} = 3.$
	\medskip
First, let us consider $A$. Let $X = Ax \simeq A/(Ax^2+Ay+Az)$ 
and $Z = Ay \simeq A/(Ax+Ay),$
these are indecomposable modules of length 2. 
We claim that for $M\in \add\{S,X,Z\},$ {\it we have $\Omega M \in \add\{S,X,Z\}$
and $t_2(M) = 3t_1(M)-t_0(M).$}
	\smallskip
Proof. We can assume that $M$ is indecomposable, thus $M$ is one of $S,X,Z.$ 
We have $\Omega S = X\oplus Z \oplus S;$ second, we have 
$\Omega X = Ax^2\oplus Z \oplus Az \simeq
Z\oplus S^2,$ and finally 
$\Omega Z = X\oplus Z$. This shows already that $\Omega M \in\add\{S,X,Z\}$.
It follows that $\Omega^2 S = \Omega(X\oplus Z \oplus S) \simeq X^2\oplus Z^3\oplus S^3$,
therefore $t_2(S) = t(\Omega^2 S) = 8$. Since $t_1(S) = 3$ and $t_0(S) = 1$, we have
$t_2(S) = 3t_1(S)-t_0(S).$ Next, $\Omega^2 X = \Omega(Z\oplus S^2) 
\simeq X^3\oplus Z^3\oplus S^2,$ 
therefore $t_2(X) = t(\Omega^2 X) = 8$. Since $t_1(X) = 3$ and $t_0(X) = 1$, we have
$t_2(X) = 3t_1(X)-t_0(X).$ Finally, $\Omega^2 Z = \Omega(X\oplus Z) 
\simeq X^2\oplus Z^2\oplus S,$ 
therefore $t_2(Z) = t(\Omega^2 Z) = 5$. Since $t_1(X) = 2$ and $t_0(X) = 1$, we have
$t_2(Z) = 3t_1(Z)-t_0(Z).$ 
$\s$
	\smallskip
By induction, we see that $\Omega^n(S) \in \add\{S,X,Z\}$ and that the numbers
$b_n = t_n(S)$ satisfy the recursion $b_{n+2} = 3b_{n+1}-b_n$ for all $n\ge 0$.
Since $b_0 = 1$ and $b_n = 3$, it follows that the numbers $b_n$ are the even-index
Fibonacci numbers $1,\ 3,\ 8,\ 21,\ 55,\ 144,\ \dots$ and therefore 
$\gamma_A = \gamma(S) = \limsup_n b_n =
\phi^2 = \psi$, where $\phi = \frac12(1+\sqrt 5)$ is the golden ratio.
	\medskip
Second, we consider $A'$. Let $X = A'x \simeq A/(Ax^2+Ay+Az)$ and $Y = A'y \simeq
A/(Ax+Ay^2+Az),$ these are indecomposable modules of length 2. 
We claim that {\it for $M\in \add\{S,X,Y\},$ we have $\Omega M \in \add\{S,X,Y\}$
and $t_1(M) = 3t_0(M)$.}
	\smallskip
Proof. We can assume that $M$ is indecomposable, thus $M$
is one of the modules $S,X,Y$. We have $\Omega S = J = X\oplus Y \oplus S$, 
and $|\top \Omega S| = 3.$ We have $\Omega X = Ax^2 \oplus Y \oplus Az \simeq Y\oplus S^2,$
thus $|\top \Omega X| = 3.$ And similarly, $\Omega Y \simeq X\oplus S^2,$ and
thus $|\top \Omega Y| = 3.$ $\s$
	\smallskip
By induction, $\Omega^n S$ belongs to $\add\{S,X,Y\}$
and $t_n(S) = |\top\Omega^n(S)| = 3^n.$ Therefore $\gamma_{A'} = \gamma(S) = 3.$
	\bigskip\bigskip 
\vfill\eject
%%%%%%%%%%%%%%%%%%%%%%%%%%%%%%%%%%%%%%%%%%%%%%%%%%%%%%%%%%%%%%%%%%%
%%%%%%%%%%%%%%%%%%%%%%%%%%%%%%%%%%%%%%%%%%%%%%%%%%%%%%%
{\bf References.}
	\medskip
\item{[AE]}  J. L. Alperin, L. Evens. Representations, resolutions,
  and Quillen's dimension theorem. J. Pure Appl. Algebra 22 (1981), 1--9.
\item{[AIS]} L. L. Avramov, S. B. Iyenga r, L. M. \c Sega. Free resolutions
  over short local rings. J. London Math. Soc. 78 (2008), 459--476.
\item{[BH]} W. Bruns, J. Herzog. Cohen-Macaulay Rings. Cambridge Studies 
   in advanced mathematics 39. Cambridge University Press (1993).
\item{[HI]} J. Herzog, S. Iyengar. Koszul modules. J. Pure Appl. Algebra
  201 (2005), 154--188.
\item{[L]} J. Lescot. Asymptotic properties of Betti numbers of modules over
  certain rings. J. Pure Appl. Algebra 38 (1985), 287--298.
\item{[RZ]}  C. M. Ringel, P. Zhang. Gorenstein-projective modules over short local
 algebras. To appear.  arXiv:1912.02081.
	\bigskip\bigskip
%%%%%%%%%%%%%%%%%%%%%%%%%%%%%%%%%%%%%%%%%%%%%%%%%%%%%%%
{\baselineskip=1pt
\rmk
C. M. Ringel\par
Fakult\"at f\"ur Mathematik, Universit\"at Bielefeld \par
POBox 100131, D-33501 Bielefeld, Germany  \par
ringel\@math.uni-bielefeld.de
\medskip

P. Zhang \par
School of Mathematical Sciences, Shanghai Jiao Tong University \par
Shanghai 200240, P. R. China.\par
pzhang\@sjtu.edu.cn\par}

\bye